\documentclass{amsart}

\usepackage{amscd,amssymb,amsmath,graphicx,verbatim}
\usepackage[dvips]{hyperref}
\usepackage{textcomp}
\usepackage[OT1,T1]{fontenc}

\theoremstyle{definition}

\theoremstyle{remark}

\numberwithin{equation}{section}

\begin{document}

\title[RH: a special case of the Riesz and Hardy-Littlewood
wave]{Riemann Hypothesis: a special case of the Riesz and Hardy-Littlewood
wave and a numerical treatment of the Baez-Duarte coefficients up
to some billions in the k-variable}
\author{Stefano Beltraminelli}
\address{S. Beltraminelli, CERFIM, Research Center for Mathematics
and Physics, PO Box 1132, 6600 Locarno, Switzerland}
\email{stefano.beltraminelli@ti.ch}
\author{Danilo Merlini}
\address{D. Merlini, CERFIM, Research Center for Mathematics and
Physics, PO Box 1132, 6600 Locarno, Switzerland}
\email{merlini@cerfim.ch}
\label{I1}
\date{17 September 2006}
\begin{abstract}
We consider the Riesz and Hardy-Littlewood wave i.e. a ``critical
function'' whose behaviour is concerned with\ \ the possible truth
of the Riemann Hypothesis (RH). The function is studied numerically
for the case $\alpha =\frac{15}{2}$ and $\beta =4$ in some range
of the critical strip, using Maple 10.\ \ \ \ 

In the experiments, $N=2000$ is the maximum argument used in the
M\"obius function appearing in $c_{k}$ i.e. the coefficients of
Baez-Duarte, in the representation of the inverse of the Zeta function
by means of the Pochammer's polynomials.

The numerical results give some evidence that the critical function
is bounded for $\mathfrak{R}( s) >\frac{1}{2}$ and such an ``evidence''
is stronger in the region $\mathfrak{R}( s) >\frac{3}{4}$ where
the wave seems to decay slowly. This give further support in favour
of the absence of zeros of the Riemann Zeta function in some regions
of the critical strip ($\mathfrak{R}( s) >\frac{3}{4}$) and a (weaker)
support in the direction to believe that the RH may be true ($\mathfrak{R}(
s) >\frac{1}{2}$). 
\end{abstract}
\maketitle

\section{Introduction}

The starting point of this note is the representation of the reciprocal
of the Riemann Zeta function by means of the Pochammer's polynomials
$P_{k}( z) $ (where {\itshape z} is a complex variable), whose coefficients
$c_{k}$ have been introduced by Baez-Duarte for the Riesz case ($\alpha
=\beta =2$). For the study of the coefficients $c_{k}$, some recent
analytical as well as numerical results have been obtained [2, 3, 4, 5, 6, 7, 8, 9].
For a rigorous treatment of the M\"untz formula to the finding of
new zero free regions of the Riemann Zeta function, the reader may
consult the work of Albeverio and Cebulla [1].

Using the Baez-Duarte approach, the representation of $\frac{1}{\zeta
( s) }$ may be obtained for a family of a two parameter Pochammer's
polynomials (parameters $\alpha$ and $\beta$ [4]) and reads:
\begin{equation}
\frac{1}{\zeta ( s) }=\sum \limits_{k=0}^{\infty }c_{k}( \alpha
,\beta ) P_{k}( s,\alpha ,\beta ) ,
\end{equation}

\noindent where
\begin{gather}
\begin{array}{rl}
 P_{k}( s,\alpha ,\beta )  & :=\prod \limits_{r=1}^{k}\left( 1-\frac{\frac{s-\alpha
}{\beta }+1}{r}\right) 
\end{array}
\end{gather}
\begin{equation}
\begin{array}{rl}
 c_{k}( \alpha ,\beta )  & :=\sum \limits_{n=1}^{\infty }\frac{\mu
( n) }{n^{\alpha }}{\left( 1-\frac{1}{n^{\beta }}\right) }^{k}.
\end{array}
\end{equation}

The expression for $c_{k}$ we will use in our computation is given
by:
\begin{equation}
c_{k}( \alpha ,\beta ) =\sum \limits_{n=1}^{N}\frac{\mu ( n) }{n^{\alpha
}}e^{-\frac{k}{n^{\beta }}},
\end{equation}

\noindent which for large {\itshape k} is a correct formula for
the $c_{k}$ given above (see Appendix). From a theorem of Baez-Duarte
[2, 3], an important inequality concerning the Pochammer's polynomials
is given by:
\begin{equation}
\left| P_{k}( z) \right| \leq A k^{-\mathfrak{R}( z) }.
\end{equation}

The inequality, when applied to our family of Pochammer's polynomials
gives:
\begin{equation}
\left| P_{k}( \frac{s-\alpha }{\beta }+1) \right| \leq A k^{-\left(
\frac{\mathfrak{R}( s) -\alpha }{\beta }+1\right) }.
\end{equation}

From this it follows [2, 3, 4] that the RH will be true, i.e.
that $\frac{1}{\zeta ( s) }$ in the representation above will be
different from infinity (no zero of $\zeta ( s) $ for $\mathfrak{R}(
s) >\rho $) if $c_{k}k^{\frac{\alpha -\rho }{\beta }}\leq C$. For
the numerical study it is convenient to introduce the variable $x=\log
( k) $, in term of which we define the critical function corresponding
to $\alpha$ and $\beta$. This is given by:
\begin{equation}
\psi ( x;\alpha ,\beta ,\rho ) :=e^{\frac{\alpha -\rho }{\beta }x}\sum
\limits_{n=1}^{2000}\frac{\mu ( n) }{n^{\alpha }}e^{-\frac{e^{x}}{n^{\beta
}}}.
\end{equation}

2000 is the maximum argument {\itshape N} used in these experiments,
which for the special case we treat ($\alpha =\frac{15}{2}$ and
$\beta =4$) $\psi$ will be calculated up to $x=30$ (this corresponds
to $k=e^{30}=1.06865\times {10}^{13}$).

Before we present the results of our numerical experiments for various
values of $\rho$ (for $\rho =1,\frac{7}{8},\frac{3}{4},\frac{5}{8},\frac{1}{2},\frac{3}{8},\frac{3}{10}$)
it is important to give the explicit expression of the contribution
of the non trivial ($\psi _{\mathit{nt}}$) and also of the trivial
zeros ($\psi _{t}$) to the critical function defined above for the
general case $\alpha$ and $\beta$, following the expression given
by Baez-Duarte for the case $\alpha =\beta =2$ [2]. For the non
trivial zeros, in the variable $x=\log ( k) $, at $\rho$, it is
given by:
\begin{equation}
\psi _{\mathit{nt}}( x;\alpha ,\beta ,\rho ) =\frac{1}{\beta }\sum
\limits_{z}^{ }\frac{e^{\frac{i \mathfrak{I}( z) }{\beta }x}\Gamma
( -\frac{\mathfrak{R}( z) +i \mathfrak{I}( z) -\alpha }{\beta })
}{\zeta ^{\prime }( z) },
\end{equation}

\noindent where {\itshape z} is any nontrivial zero of $\zeta (
s) $. In our experiments we will limit to the contribution of the
first two lower zeros given experimentally by $z_{1}=\frac{1}{2}+i
14.134725...$ and $z_{2}=\frac{1}{2}+i 21.022040...$ and the complex
conjugate of them. The corresponding contribution will be denoted
by $r_{1}( x) $ (from $z_{1}$ and ${\overline{z}}_{1}$) and $r_{2}(
x) $ (from $z_{2}$ and ${\overline{z}}_{2}$).

The contribution of the trivial zeros $z=-2n$ for every integers
{\itshape n}, to the critical function is given by:
\begin{equation}
\psi _{\mathit{t}}( x;\alpha ,\beta ,\rho ) =\frac{1}{\beta }\sum
\limits_{n=1}^{\infty }\frac{e^{-\frac{2n+\rho }{\beta }x}\Gamma
( \frac{\alpha +2n}{\beta }) }{\zeta ^{\prime }( -2n) },
\end{equation}

\noindent where a summation until $N=20$ will be sufficient.

So, in our calculations we will set $\alpha =\frac{15}{2}$ and $\beta
=4$ in the above formulas, for any value of $\rho$ we shall consider.
The contribution $\psi _{t}$ for $\rho $ will be indicated with
$g_{\rho }( x) $. Below we present the results of our numerical
experiments performed using Maple 10, where as anticipated the maximum
argument in the M\"obius function present in the definition of the
critical function (essentially $c_{k}$), is $N=2000$. The fluctuations
errors around 2000 will be specified in the Appendix.

\section{Numerical experiments}

In Fig. 1 we give the plot of the two functions $\psi ( x;\frac{15}{2},4,\frac{1}{2})
-r_{1}( x) -r_{2}( x) $ and $g_{1/2}( x) $ up to $x=30$ which shows
a very good agreement. Notice that we have taken into account only
the contribution of the first two nontrivial zeros in the Baez-Duarte
asymptotic formula for the $c_{k}$. For the Riesz case ($\alpha
=\beta =2$), the contribution of the trivial zeros to the $c_{k}$
have been treated by Maslanka using the Rice's integrals [8].
\begin{figure}[h]
\begin{center}
\includegraphics[width=8cm]{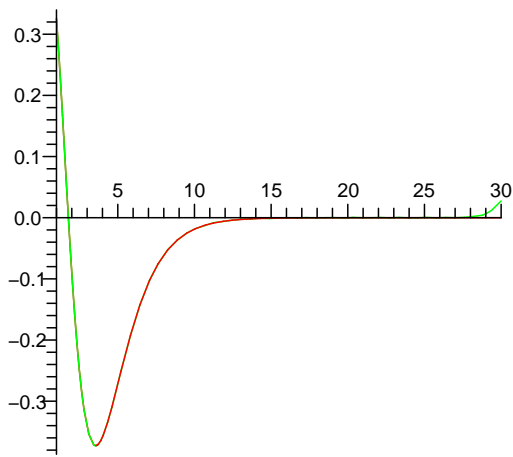}

\end{center}
\caption{Plot of $\psi ( x) -r_{1}( x) -r_{2}( x) $ [red] and $g_{1/2}(
x) $ [green] up to $x=30$}

\end{figure}

In Fig. 2 we present the the plot of the two functions $\psi ( x;\frac{15}{2},4,\frac{1}{2})
-g_{1/2}( x) $ and $r_{1}( x) +r_{2}( x) $ up to $x=30$ which shows
not only a good agreement but also the oscillatory behaviour of
the contribution of the first two nontrivial zeros. 
\begin{figure}[h]
\begin{center}
\includegraphics[width=8cm]{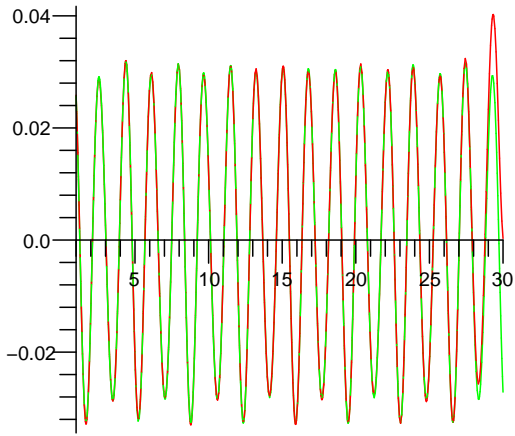}

\end{center}
\caption{Plot of $\psi ( x) -g_{1/2}( x) $ [red] and $r_{1}( x)
+r_{2}( x) $ [green] up to $x=30$}

\end{figure}

In the next Fig. 3 we present the plots of some critical functions
($\psi _{\rho }$) corresponding to different values of $\rho$ using
(1.7) and this without any comparison with the Baez-Duarte asymptotic
expansion considered above. It is to be noted that all functions
$\psi _{\rho }$ has the same zeros and we observe that there is
a well marked evidence that for $\rho >\frac{1}{2}$ increasing to
1 the amplitudes decay while for $\rho <\frac{1}{2}$ the amplitudes
grow. These functions have been indicated with $\psi _{1}$, $\psi
_{7/8}$, $\psi _{3/4}$, $\psi _{5/8}$, $\psi _{1/2}$, $\psi _{3/8}$,
$\psi _{3/10}$ respectively.

It should be said that $\psi _{3/8}$ and $\psi _{3/10}$, we have
considered, have no relation with the representation of $\frac{1}{\zeta
( s) }$ which is valid only for $\mathfrak{R}( s) >\frac{1}{2}$.
The two functions help only to visualize that $\psi _{1/2}$ is the
borderline for the critical functions decaying for $\mathfrak{R}(
s) >\frac{1}{2}$ as suggested by our numerical experiments up to
$x=30$. It should also be added that from the duality relation (Riemann
symmetry of the Zeta function), given by:\ \ \ \ \ 
\begin{equation}
\frac{1}{\zeta ( 1-s) }=\pi ^{s-\frac{1}{2}}\frac{\Gamma ( \frac{1-s}{2})
}{\Gamma ( \frac{s}{2}) }\frac{1}{\zeta ( s) },
\end{equation}

\noindent it follows that the right hand side of (2.1) ensures a
representation of $\frac{1}{\zeta ( s) }$ via the Pochammer's polynomials
in the region $0<\mathfrak{R}( s) <\frac{1}{2}$. 
\begin{figure}[h]
\begin{center}
\includegraphics[width=8cm]{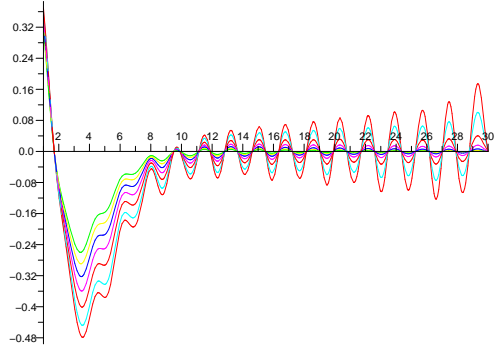}

\end{center}
\caption{Plot of $\psi _{\rho }$ for $\rho =1,\frac{7}{8},\frac{3}{4},\frac{5}{8},\frac{1}{2},\frac{3}{8},\frac{3}{10}$
up to $x=30$, in order of increasing amplitudes}

\end{figure}

In Fig. 4 we present the result for a special case where we allow
a slower decrease in the critical function (see addendum in the
exponent of the critical function), which is the same as to say
that we ask only for a slower decay of $c_{k}$, at $\rho =\frac{1}{2}$
i.e. of the type $c_{k}=\frac{A \log ( k) }{k^{\frac{7}{4}}}$ for
the case considered. This is not the same as to ask that RH is true
or that RH is true with nontrivial zeros which are simple [3].
It is a case in between the two.

In this case the critical function (indicated with $\psi _{1/2+}$)
is explicitly given by:
\begin{equation}
\psi _{1/2+}( x) =e^{\frac{7}{4}x-\log ( x) }\sum \limits_{n=1}^{2000}\frac{\mu
( n) }{n^{\frac{15}{2}}}e^{-\frac{e^{x}}{n^{4}}}.
\end{equation}
\begin{figure}[h]
\begin{center}
\includegraphics[width=8cm]{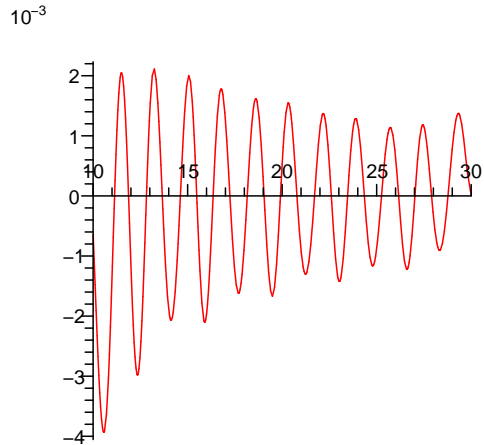}

\end{center}
\caption{Plot of $\psi _{1/2+}$ }

\end{figure}

Here there is more evidence that the amplitude of the wave at $\rho
=\frac{1}{2}$ is decreasing with $x=\log ( k) $. The experiments
of Fig. 3 give in any cases a stronger evidence: that for $\rho
>\frac{3}{4}$ the amplitudes of the waves are decaying, and thus
are bounded in amplitude by a constant. This is a symptom of the
absence of nontrivial zeros in the critical segment $\frac{3}{4}<\rho
<1$.

In the last experiment we set $\rho =\frac{3}{4}$ and compare $\psi
_{3/4}$ with the asymptotic expression of Baez-Duarte: for the trivial
zeros we set $\rho =\frac{3}{4}$ in the above formula, for the nontrivial
zeros (the two we consider) we keep the same value of $\mathfrak{I}(
z) $ but we assume that their real part is $\mathfrak{R}( z) =\frac{3}{4}$.
The plot in Fig. 5 of the function $\psi _{3/4}( x) $ and of $g_{3/4}(
x) +r_{1}( x) +r_{2}( x) $ are clearly different: in $\psi _{3/4}$
there is the trace via the M\"obius function of where the nontrivial
zeros are located and thus the amplitude is decaying. In the second
function, the two considered zeros are supposed to have $\mathfrak{R}(
z) =\frac{3}{4}$ and the wave which appears seems to have a constant
amplitude as in the case $\psi _{1/2}$ which of course would be
sufficient to ensure the truth of the RH.
\begin{figure}[h]
\begin{center}
\includegraphics[width=8cm]{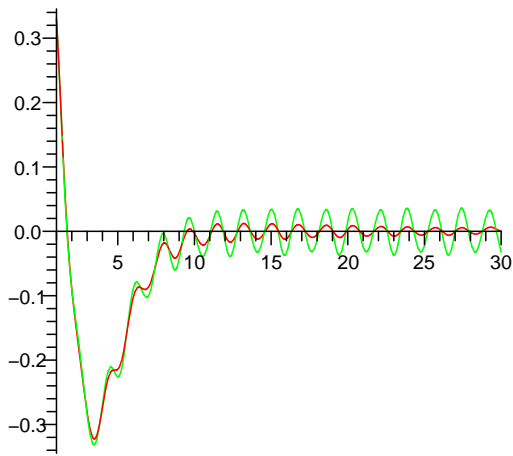}

\end{center}
\caption{Plots of the functions $\psi _{3/4}( x) $ [red] and $g_{3/4}(
x) +r_{1}( x) +r_{2}( x) $ [green]}

\end{figure}

In Appendix we analyse a (weak) stability property of our results
obtained with $N=2000$ in the M\"obius function and give some indications
why the waves for $\rho =\frac{3}{4}$ should be decaying, thus ensuring
more credibility on the absence of zeros of the Riemann Zeta function
in the segment $\frac{3}{4}<\rho <1$.

\section{Conclusions}

In this work we have analyzed numerically the behaviour of the Riesz
and Hardy-Littlewood wave (the critical function) in some details
for the case $\alpha =\frac{15}{2}$ and $\beta =4$ in the region
up to about $k=10'000$ milliard for various values of $\rho$ in
the critical segment $\frac{1}{2}<\rho <1$. In the variable $x=\log
( k) $ up to 14 oscillations have been detected whose amplitude
has been compared with the one calculated with the expansion of
Baez-Duarte using the trivial zeros and only the two lower nontrivial
zeros. The agreement is satisfactory and the results give some indication
in the direction to believe that at least for $\mathfrak{R}( s)
=\rho >\frac{3}{4}$ there are no nontrivial zeros of $\zeta ( s)
$ since in the representation $\frac{1}{\zeta ( s) }$ seems to remain
bounded. In addition, we have given some evidence that a slow decay
of $c_{k}$ like $\frac{\log ( k) }{k^{\frac{7}{4}}}$ in between
to the decay $\frac{k^{\epsilon }}{k^{\frac{7}{4}}}$ (RH for the
model) and $\frac{A}{k^{\frac{7}{4}}}$ (RH for the model with simple
zeros) is possible in the range $\log ( k) <30$. A further study
in this direction but by means of two new representations of the
Zeta function with coefficients $b_{k}$ and $d_{k}$ and other type
of oscillations will be presented in the near future.

\appendix

\section{}

We consider the critical function $\psi _{3/4}( x) $ obtained with
$N=2000$ (maximum argument in the M\"obius function appearing in
the Baez-Duarte definition of the $c_{k}$). We will suppose that
the numerical results are given with good accuracy; we now ask:
if we increase {\itshape N} from 2000 up to ${10}^{6}$ in a ideal
experiment, what will be the change of the critical function in
the range $x<30$?
\begin{gather*}
\psi _{3/4}( x;N=2000) = e ^{\frac{27}{16}x}\sum \limits_{n=1}^{2000}\frac{\mu
( n) }{n^{\frac{15}{2}}} e ^{-\frac{ e ^{x}}{n^{4}}}
\end{gather*}
\[
\psi _{3/4}( x;N={10}^{6}) = e ^{\frac{27}{16}x}\sum \limits_{n=1}^{{10}^{6}}\frac{\mu
( n) }{n^{\frac{15}{2}}} e ^{-\frac{ e ^{x}}{n^{4}}}.
\]

The difference $\Delta$ between the two functions is bounded ($|\mu
( n) |\leq 1$) by:
\[
\Delta \leq  e ^{\frac{27}{16}x}\sum \limits_{n=2000}^{{10}^{6}}\frac{1}{n^{\frac{15}{2}}}
e ^{-\frac{ e ^{x}}{{10}^{24}}}.
\]

If we ask that $\Delta$ will be smaller then say ${10}^{-6}$ time
0.015 which is about the value of the amplitude of the wave in the
range $x\leq 30$, obtained with $N=2000$, we have:
\[
\Delta \leq  e ^{\frac{27}{16}x} e ^{-\frac{ e ^{x}}{{10}^{24}}}(
\zeta ( \frac{15}{2}) -\zeta ( \frac{15}{2};N=2000) ) \leq 0.015\cdot
{10}^{-6}.
\]

The difference between the Zetas is estimated to:
\[
\operatorname*{\int }\limits_{2000}^{\infty }\frac{1}{x^{\frac{15}{2}}}dx=\frac{2}{13}{2000}^{-\frac{13}{2}}=\frac{2}{65}{10}^{-26}.
\]

And the inequality takes the form:
\[
\frac{27}{16}x- e ^{\frac{x}{{10}^{24}}}+\log ( \frac{2}{65}) -26
\log ( 10) +6 \log ( 10) -\log ( 0.015) \leq 0,
\]

\noindent with the solution $x\leq 27$. Thus for $x\leq 27$, the
amplitudes will change at most ${10}^{-6}$ time of its value 0.015.
This shows some stability in the numerical experiments as {\itshape
N} increases in a ideal experiment. Of course this is independent
of how many zeros are employed in the Baez-Duarte estimation.

Finally, application of the crude inequality as above shows [5]
that for larger values of {\itshape x}, with {\itshape N} sufficiently
big, the numerical results of an ideal experiment using the M\"obius
function in $c_{k}$ will give a value of the amplitude for example
smaller then $\frac{1}{2}$ of 0.015 (still with $N=2000$ at $x>35$,
with $N={10}^{9}$ at $x>88$), indicating that the value of the amplitude
becomes possibly smaller. This is an indication in the direction
to believe that for $\rho >\frac{3}{4}$ there are no nontrivial
zeros of the Zeta function.

To end up with the Appendix it should be remarked that in our experiments
we have used the formula for the $c_{k}$ [4], given by:
\[
{\hat{c}}_{k}=\sum \limits_{n=1}^{\infty }\frac{\mu ( n) }{n^{\alpha
}} e ^{-\frac{k}{n^{\beta }}},
\]

\noindent instead of the formula: 
\[
c_{k}=\sum \limits_{n=1}^{\infty }\frac{\mu ( n) }{n^{\alpha }}{\left(
1-\frac{1}{n^{\beta }}\right) }^{k}.
\]

Again, as above, the crude inequality $|\mu ( n) |\leq 1$ may be
used to show that the difference between the two functions i.e.
the fluctuations become smaller as {\itshape k} get bigger and depends
on $\alpha$ and $\beta$. In fact they behave unconditionally as:
\[
A\sum \limits_{n=1}^{\infty }\frac{1}{n^{\alpha +\beta }}{\left(
1-\frac{1}{n^{\beta }}\right) }^{k}\leq \frac{C}{k^{\frac{\alpha
+\beta -1}{\beta }}}.
\]

To see that, let $\Delta =|{\hat{c}}_{k}-c_{k}|$ then:
\[
\Delta \leq \sum \limits_{n=1}^{\infty }\frac{\left| \mu ( n) \right|
}{n^{\alpha }}\left(  e ^{-\frac{k}{n^{\beta }}}-{\left( 1-\frac{1}{n^{\beta
}}\right) }^{k}\right) \leq \sum \limits_{n=1}^{\infty }\frac{1}{n^{\alpha
}}\left(  e ^{-\frac{k}{n^{\beta }}}-{\left( 1-\frac{1}{n^{\beta
}}\right) }^{k}\right) ,
\]

\noindent since $ e ^{-\frac{k}{n^{\beta }}}\geq {(1-\frac{1}{n^{\beta
}})}^{k}$ and passing to the continuous variable {\itshape x}, the
contribution of the second integral is given by [4]: 
\[
\operatorname*{\int }\limits_{1}^{\infty }\frac{1}{x^{\alpha }}{\left(
1-\frac{1}{x^{\beta }}\right) }^{k}dx=\frac{1}{\beta }\frac{\Gamma
( \frac{\alpha -1}{\beta }) \Gamma ( k+1) }{\Gamma ( \frac{\alpha
-1}{\beta }+k+1) },
\]

\noindent while the first is given by:
\[
\operatorname*{\int }\limits_{1}^{\infty }\frac{ e ^{-\frac{k}{x^{\beta
}}}}{x^{\alpha }}dx=\frac{1}{\beta  k^{\frac{\alpha -1}{\beta }}}\Gamma
( \frac{\alpha -1}{\beta }) .
\]

At large {\itshape k} the fluctuation behaves like the difference,
i.e. as:
\[
\Delta \leq \frac{C}{k^{\frac{\alpha +\beta -1}{\beta }}}.
\]

For the model under consideration the decay is as $\frac{C}{k^{\frac{21}{8}}}$
and is stronger then in the usual Riesz case ($\alpha =\beta =2$)
where an early more detailed calculation gives a decay like $\frac{C}{k^{\frac{3}{2}}}$
[6]. 

Finally it should be added that the general upper bound for $\Delta$
is related to the discrete derivative of the Baez-Duarte coefficients
given by:
\begin{align*}
c_{k}-c_{k+1}&=\sum \limits_{n=1}^{\infty }\frac{\mu ( n) }{n^{\alpha
}}\left( {\left( 1-\frac{1}{n^{\beta }}\right) }^{k}-{\left( 1-\frac{1}{n^{\beta
}}\right) }^{k+1}\right) 
\\%
 &=\sum \limits_{n=1}^{\infty }\frac{\mu ( n) }{n^{\alpha }}{\left(
1-\frac{1}{n^{\beta }}\right) }^{k}\left( 1-1+\frac{1}{n^{\beta
}}\right) =c_{k}( \alpha +\beta ,\beta ) ,
\end{align*}

\noindent which unconditionally are bounded by $\frac{C}{k^{\frac{\alpha
+\beta -1}{\beta }}}$ as above [4].

In the same way
\[
-\frac{d}{dk}\sum \limits_{n=1}^{\infty }\frac{\mu ( n) }{n^{\alpha
}} e ^{-\frac{k}{n^{\beta }}}=\sum \limits_{n=1}^{\infty }\frac{\mu
( n) }{n^{\alpha +\beta }} e ^{-\frac{k}{n^{\beta }}},
\]

\noindent which gives the same decay since the function is equal
to $c_{k}( \alpha +\beta ,\beta ) $ as above.\ \ 

At large {\itshape k} we also have [4]: 
\[
c_{k}\approx \sum \limits_{p=0}^{\infty }\frac{c_{p}k^{p} e ^{-k}}{p!}
\]

\noindent a Poisson like distribution for the coefficients $c_{k}$.

\end{document}